\documentclass[11pt,leqno]{article}
\usepackage{amsthm,amsmath,amssymb,graphicx,amsfonts}
\usepackage{dsfont}
\RequirePackage{hypernat}
\usepackage[aboveskip=10pt,font=small,labelfont=bf]{caption} 
\usepackage{subfig}
\usepackage{float}
\usepackage[section]{placeins}
\RequirePackage[colorlinks,citecolor=blue,urlcolor=blue]{hyperref}
\newcounter{pb}
\numberwithin{equation}{section}

\newtheorem{theorem}{Theorem}[section]
\newtheorem{corollary}[theorem]{Corollary}

\newtheorem{definition}[theorem]{Definition}

\newtheorem*{KPZ}{KPZ formula}

\newtheorem{proposition}[theorem]{Proposition}

\newtheorem{remark}[theorem]{Remark}
\newtheorem*{remark*}{Remark}

\def \Cv{\mathrm{Cov}}
\def \R{ \mathbb R }

\def \one{ {\rm 1}\mkern-4.5mu{\rm I} }

\def \E{ \mathbb  E }

\def \P{ \mathbb P  }
\def \Q{ \mathbb Q  }

\def \T{\mathds T}
\begin{document}

\title{Optimal transport for multifractal random measures. Applications.}
\maketitle
\begin{center}
{R\'emi Rhodes, Vincent Vargas \\
\footnotesize

 CNRS, UMR 7534, F-75016 Paris, France \\
  Universit{\'e} Paris-Dauphine, Ceremade, F-75016 Paris, France} \\

{\footnotesize \noindent e-mail: \texttt{rhodes@ceremade.dauphine.fr},\\
\texttt{vargas@ceremade.dauphine.fr}}
\end{center}

\begin{abstract}
In this paper, we study optimal transportation problems for multifractal random measures. Since these measures are much less regular than  optimal transportation theory requires, we introduce a new notion of transportation which is intuitively some kind of multistep transportation. Applications are given for construction of multifractal random changes of times and to the existence of random metrics, the volume forms of which coincide with the multifractal random measures. This study is motivated by recent problems in the KPZ context.
\end{abstract}
\vspace{1cm}
\footnotesize

%\noindent{\bf Short Title.}

\noindent{\bf Key words or phrases:} Random measures, multifractal processes, optimal transport, metric, KPZ.


\normalsize

%%%%%%%%%%%%%%%%%%%%%%%%%%%%%%%%%%%%%%%%%%%%%%%%%%%%%%%%%%%%%%%%%%%%%%%%%%%%%%%%%%%%%%%%%%

\section{Introduction}

On the Borelian subsets of $\R^m$, consider a measure $M$ formally defined by
\begin{equation}\label{measureM}
M(A)=\int_Ae^{\gamma X(x)-\frac{\gamma^2}{2}\E[X(x)^2]}\,dx,
\end{equation} 
where $\gamma>0$ is a parameter,  $(X(x))_{x\in\R^m}$ is a Gaussian distribution with covariance function given by
\begin{equation}\label{cov}
\Cv(X(x),X(y))=\ln_+\frac{T}{|y-x|}+g(|y-x|),
\end{equation} 
and  $g$ is a continuous bounded function. Actually, whatever the function $g$ is, the part that really matters in \eqref{cov} is the logarithm part. Such measures are called log-normal multifractal random measures and were first rigorously defined in \cite{kah}.

The above situation can be generalized to the situation 
\begin{equation}\label{levyM}
M(A)= \int_Ae^{\omega(x)}\,dx
\end{equation}
where the process $\omega$ is a suitable Levy distribution: the resulting measures are called {\it log-infinitely divisible multifractal random measures}, MRM for short (see subsection \ref{cons:loglevy} for a reminder of the construction). Such measures exhibits interesting properties like stationarity, isotropy, long-range dependence, fat tail distribtions and, because of their log part in \eqref{cov} (or a suitable generalization for L\'evy distributions), possess a remarkable scaling property, the so-called  {\bf stochastic scale invariance}: \begin{equation}\label{esip}
(M(\lambda A))_{A\subset B(0,T)}\stackrel{\text{law as } \lambda\to 0}{\simeq}(\lambda^m e^{\Omega_\lambda}  M(A))_{A\subset B(0,T)},
\end{equation} 
where $\Omega_\lambda$ is an infinitely divisible random variable, independent of $ (M(A))_{A\subset B(0,T)}$. When $M$ satisfies the above relation \eqref{esip} with $=$ instead of $\simeq$, we will say that $M$ satisfies the exact stochastic scale invariance property (or $M$ is ESSI for short).

The purpose of this paper is to investigate (optimal) transportation problems associated to these measures.
A transport map between two probability measures $\mu,\nu $ is a map that pushes $\mu$ forward to $\nu$. The transport map is said to be optimal if it realizes the infimum of a cost functional among all the possible transport maps. For usual cost functionals, existence and uniqueness of an optimal transport map are strongly connected to the regularity of the measures $\mu,\nu$. Concerning MRM, their regularity is much weaker than that required in optimal transportation theory. So we give new notions of transportation that can be applied to MRM. This study is highly motivated by the construction of multifractal random changes of time and the construction of metric spaces the volume form of which is given by the MRM, as explained in subsection \ref{sec:appl}. The latter construction is related to a recent problem in the KPZ context raised by K. Izyurov at the conference "Les Diablerets".

%%%%%%%%%%%%%%%%%%%%%%%%%%%%%%%%%%%%%%%%%%%%%%%%%%%%%%%%%%%%%%%%%%%%%%%%%%%%%%%%%%%%%%%%%%%%%%%%%%%%%%%%%%%%%%%%%%%%%%%%%%%%%%%%%%%%%%%%%%%%%%%%%%%%%%%%%%%%%%%%%%%%
\section{Background and main results}
%%%%%%%%%%%%%%%%%%%%%%%%%%%%%%%%%%%%%%%%%%%%%%%%%%%%%%%%%%%%%%%%%%%%%%%%%%%%%%%%%%%%%%%%%%%%%%%%%%%%%%%%%%%%%%%%%%%%%%%%%%%%%%%%%%%%%%%%%%%%%%%%%%%%%%%%%%%%%%%%%%%%
In this section we first give the basic background in order to state rigorously the main results of the paper. Since the function $g$ in \eqref{cov} does not play a part in what follows, we focus on the case where the measure $M$ satifies the exact scale invariance property.

\subsection{Reminder of the construction of ESSI MRM}\label{cons:loglevy}
%%%%%%%%%%%%%%%%%%%%%%%%%%%%%%%%%%%%%%%%%
We present below the generalization of \eqref{measureM} to the situation where $X$ is a L\'evy distribution. For further details, the reader is referred to \cite{rv2}. To characterize such a L\'evy distribution, we consider the characteristic function of an infinitely divisible random variable $Z$, which can be written as $\E[e^{iqZ}]=e^{\varphi(q)}$ where (L\'evy-Khintchine's formula)
$$ \varphi(q)=i bq-\frac{1}{2}\sigma^2q^2+\int_{\R^*} (e^{iqx} -1 -iq\sin(x))\,\nu(dx)$$
and $\nu(dx)$ is a so-called L\'evy measure satisfying $\int_{\R^*}\min(1,x^2)\,\nu(dx)<+\infty$. We also introduce the Laplace exponent $\psi$ of $Z$ by $\psi(q)=\varphi(-iq)$ for each $q$ such that both terms of the equality make sense, and we assume that $\psi(1)=0$ (renormalization condition), $\psi(2)<+\infty$ and $\int_{[-1,1]}|x|\,\nu(dx)<+\infty$ (sufficient conditions for existence of MRM).

Now we define the process $\omega$ of \eqref{levyM}. We remind that this process has to be stationary, isotrop  and  suitably scale invariant.  To ensure isotropy, we introduce the unitary group $G$ of $\R^m$, that is
$$G=\{M\in M_m(\R);MM^t=\textrm{I}\},$$ and $H$ its  unique right translation invariant Haar measure  with mass 1 defined on the Borel $\sigma$-algebra $\mathcal{B}(G)$.  For stationarity and scale invariance, we consider the half-space
$$ S=\{(t,y);t\in\R,y\in\R^*_+\},$$ with which we associate the measure (on the Borel $\sigma$-algebra $\mathcal{B}(S)$)
$$\theta(dt,dy)=y^{-2}dt\,dy. $$
Then we consider an independently scattered infinitely divisible random measure $\mu$ associated to $(\varphi,H\otimes \theta)$ and distributed on $G\times S$. In particular, we have $\E(e^{iq\mu(A)})=e^{\varphi(q) H\otimes\theta(A)}$.

%\begin{definition}{\bf Filtration $\mathbf{\mathcal{F}_l} $.}
%Let $\Omega$ be the probability space on which $\mu$ is defined. $\mathcal{F}_l$ is defined as the $\sigma$-algebra generated by $\{\mu(A\times B); A \subset G,B\subset S, {\rm dist}(B,\R^2\setminus S)\geq l\}$.
%\end{definition}
We combine the previous ingredients in a geometrical way to define $\omega$.
Given $T$, let us  define the function $f:\R_+\to \R$ by
$$ f(l)=\left\{
\begin{array}{ll}
 l, &  \text{ if } l\leq T   \\
 T &  \text{ if } l\geq T
\end{array}\right.$$ The cone-like subset $A_l(t)$ of $S$ is defined by
$$ A_l(t)=\{(s,y)\in S; y\geq l, -f(y)/2\leq s-t\leq f(y)/2\}.$$
For any $x\in\R^m$ and $g\in G$, we denote by $x^g_1$ the first coordinate of the vector $gx$. The cone product $C_l(x)$ is then defined as 
$$C_l(x)=\{(g,t,y)\in G\times S;(t,y)\in A_l(x^g_1)\},$$
and the process $\omega_l$ ($0<l<T$)  by $ \omega_l(x)=\mu(C_l(x))$ for $x\in \R^m$.
%For forthcoming computations, we stress that for $s,t$ real we have: 
%\begin{equation*}
% \theta(A_l(s) \cap A_l(t))=g_l(|t-s|) 
%\end{equation*}
%where $g_l:\R_+\to \R$ is given by (with the notation $x^+=\text{max}(x,0)$):
%\begin{equation*} 
%g_l(x)=\left\{
%\begin{array}{ll}
% \ln(T/l)+1-\frac{x}{l}, &  \text{ if } x \leq l   \\
% \ln^{+}(T/x) &  \text{ if } x \geq l
%\end{array}\right.
%\end{equation*}

The Radon measure $M$ is then defined as the almost sure limit (in the sense of weak convergence of Radon measures) by 
$$M(A)=\lim_{l\to 0^+} M_l(A)=\lim_{l\to 0^+} \int_Ae^{\omega_l(x)}\,dx$$
 for any Lebesgue measurable subset $A\subset \R^m$. The convergence is ensured by the fact that the family $(M_l(A))_{l>0}$ is a right-continuous positive martingale. The structure exponent of $M$ is defined by
$$\forall q \geq 0, \quad \zeta(q)= dq-\psi(q) $$ for all $q$ such that the right-hand side makes sense.
The measure $M$ is different from $0$ if and only if there exists $\epsilon>0$ such that $ \zeta(1+\epsilon)>m$, (or equivalently $\psi'(1)<m$). In that case, we have:

\begin{theorem}\label{prop:existMRM}
The measure $M$ is stationary, isotropic and satisfies the  {\bf exact stochastic scale invariance property}: for any $\lambda\in ]0,1]$, 
$$ (M(\lambda A))_{A\subset B(0,T)}\stackrel{{\rm law}}{=}(\lambda^m e^{\Omega_\lambda}  M(A))_{A\subset B(0,T)},$$ where $\Omega_\lambda$ is an infinitely divisible random variable, independent of $ (M(A))_{A\subset B(0,T)}$, the law of which is characterized by:
$$ \E[e^{iq \Omega_\lambda}]=\lambda^{-\varphi(q)}.$$
\end{theorem}
 
Furthermore, the support of such measures is full in the sense that 

\begin{proposition}\label{support}
If the measure $M$ is non-degenerate, that is $\psi'(1)<m $, we have ${\rm Supp} (M)=\R^m$. Consequently, every Borelian subset of $\R^m$ with null $M$-measure has its complement dense in $\R^m$ for the Euclidian distance.
\end{proposition} 
\noindent {\bf Proof.}
%Since  $\overline{M}\big({\rm Supp} (\Gamma)\big)=1$,  it suffices to prove that, $\P$ a.s., the support of $\overline{M}$ is $B_R$. 
For a given ball $B(x,r)\subset \R^m$, the event $\{M(B(x,r))>0\} $ is measurable with respect to the asymptotic sigma algebra generated by $(\omega_l)_{l>0}$. By the $0$-$1$ law, it has probability $0$ or $1$. Because of the uniform integrability of the martingale $(M_l(A))_{l>0}$ for each Borelian subset $A$ with finite Lebesgue measure (denoted by $\lambda(A)$), we have $\E[M(B(x,r))]=\lambda(B(x,r))>0$.  We deduce $\P(M(B(x,r))>0)=1$. Hence,  $\P$ a.s., $M(B(x,r))>0$ for all the balls $B(x,r)$ with rational centers and radii.\qed

\begin{figure}
\centering 
%%----start of first subfigure---- 
\subfloat[$\gamma^2\simeq 0.1$]{
\label{fig:stacksub:a} %% label for first subfigure 
\includegraphics[width=6.3cm,height=6.3cm]{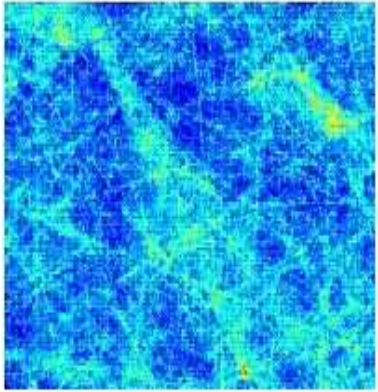}} 
\hspace{0.1\linewidth} 
%----start of second subfigure---- 
\subfloat[$\gamma^2\simeq 1.5$]{ 
\label{fig:stacksub:b} %% label for second subfigure 
\includegraphics[width=6.3cm,height=6.3cm]{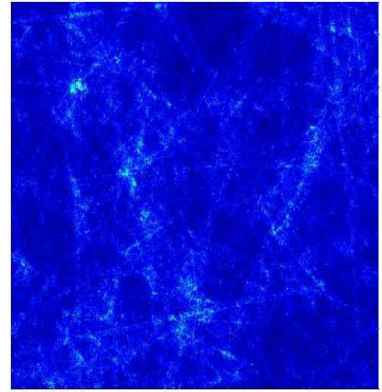}}\\[20pt] 
%%----start of third subfigure---- 
\subfloat[$\gamma^2\simeq 2$]{ 
\label{fig:stacksub:c} %% label for third subfigure 
\includegraphics[width=6.3cm,height=6.3cm]{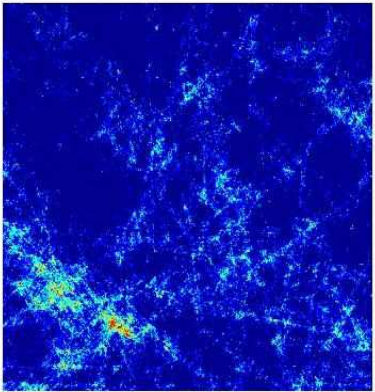}} 
\hspace{0.1\linewidth} 
%%----start of fourth subfigure---- 
\subfloat[$\gamma^2\simeq 3.9$]{ 
\label{fig:stacksub:d} %% label for fourth subfigure 
\includegraphics[width=6.3cm,height=6.3cm]{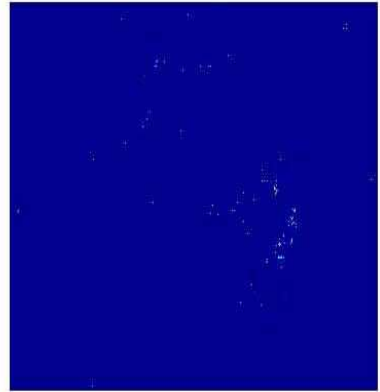}} 
\caption{Simulations of the density of a log-normal MRM with various intermittency parameters $\gamma$ appearing in \eqref{measureM}. The structure exponent matches $\xi(q)=(2+\frac{\gamma^2}{2})q-\frac{\gamma^2}{2}q^2$. In dimension $2$, the measure is non-degenerate provided that $0<\gamma^2<4$. The last two simulations are colored with a logarithmic intensity scale.}
\label{fig:stacksub} %% label for entire figure 
\end{figure}

\subsection{Optimal transport for MRM}\label{optimal}
%%%%%%%%%%%%%%%%%%%%%%%%%%%%%%

The reader may consult Appendix \ref{optimaltheory} for a brief reminder about optimal transport theory. In particular, we adopt the notations used in  Appendix \ref{optimaltheory}.
We  denote by $B_R$ the closed ball of $\R^m$ centered at $0$ with radius $R$ and by $C_R$ its Lebesgue measure. We define $\lambda_R$ as the renormalized Lebesgue measure on $B_R$ with mass $1$ and  the renormalized probability measure on $B_R$ 
$$  \overline{M}(dx)=\frac{1}{M(B_R)}M(dx).$$ 

Originally, our problem was to construct a nice (as much as possible) mapping that pushes the MRM $M$ forward to the Lebesgue measure. Of course, there are several possibilities to achieve that construction. 
Among all the possibilities, we focus on optimal transportation theory because it features many interesting qualities: measurability with respect to the randomness, isotropy, regularity,... Basic results ensure that there exists an  optimal transport map pushing $\overline{M}$ forward to the Lebesgue measure $\lambda_R$ (for a quadratic cost function) provided that the measure $\overline{M}$ does not give mass to small sets. Though  that condition is usually not satisfied by MRMs, this is true when the intermittency parameter $\psi'(1)$ is small:

%If we denote by $\pi_n$ the probability measure $({\rm Id},\chi_n)_\#\lambda_R$. This probability measure on $B_R\times B_R$ weakly converges as $n\to +\infty$ towards an optimal transference plan $\pi $ between  $\lambda_R$ and $\overline{M}$. Roughly speaking, one may convince oneself that the curve $(t,\gamma_n^{x,y}(t))_{0\leq t \leq 1}$ should converge towards a possibly multivariate curve lying in the support of $\pi$. In case $\pi$ is deterministic, that support is the graph of the gradient of some convex function $\Gamma$ so that the curve $(t,\gamma_n^{x,y}(t))_{0\leq t \leq 1}$ should converge towards the curve $$\Big(t,\chi\big((1-t)\Gamma(x)+t\Gamma(y)\big)\Big)_{0\leq t \leq 1}.$$ Hence there is much to learn of the convergence of geodesics from the existence of a unique transport map between $\overline{M}$ and $\lambda_R$. This is the purpose of our paper. Because of the singularity of the measures $\overline{M}$ and $\lambda_R$, the function $\Gamma$ has very bad regularity properties. Hence our purpose it to give a sufficient condition on the parameter $ \gamma$ of the measure $\overline{M}$ in order to have a unique optimal transport map between $\overline{M}$ and $\lambda_R$ and to give a meaning to the convergence of the curves.

\begin{theorem}\label{gamsmallevy}
When $\psi'(1)<1$, the measure $\overline{M}$ does not give mass to small sets. Hence there is a unique optimal transport map (for the Euclidian quadratic cost) that pushes the renormalized probability measure $\overline{M}$ forward to the renormalized Lebesgue measure.
\end{theorem}

We can thus apply the classical transportation theory (see Theorem \ref{th:transport}). There exist two optimal transport maps 
$$\chi:{\rm Supp}(\chi) \to{\rm Supp}(\Gamma)\quad \text{ and }\quad \Gamma:{\rm Supp}(\Gamma)\to{\rm Supp}(\chi),$$
(the supports of which  are Borelian and contained in $B_R$) that respectively push the Lebesgue measure $\lambda_R$ forward to $\overline{M}$ and vice-versa, meaning
$$\chi_\#\lambda_R =\overline{M}\quad \text{ and }\quad \Gamma_\#\overline{M}=\lambda_R.$$
They are both unique, gradients of  convex functions and bijections, inverse from each other. They are respectively $\lambda_R$ and $\overline{M}$-almost surely defined, meaning $\overline{M}\big({\rm Supp} (\Gamma)\big)=1$ and  $\lambda_R\big({\rm Supp} (\chi)\big)=1$, so that both supports ${\rm Supp} (\Gamma)$, ${\rm Supp} (\chi)$ are dense in $B_R$ for the Euclidian distance.

\vspace{2mm}
We equip $B_R$ with the Riemannian metric 
$$\forall x\in B_R,\forall u,v \in \R^m, g^e_x(u,v)=\frac{M(B_R)^2}{C_R^2}(u,v)$$ where $(\cdot,\cdot)$ denotes the usual inner product on $\R^m$.
In that way, $(B_R,g^e)$ is a Riemannian space in which the volume form matches $M(B_R)\times \lambda_R(dx)$ and the geodesic distance $d_e$ is given by the Euclidian distance on $B_R$ up to a multiplicative constant:
$$d_e(x,y)=\frac{M(B_R)}{C_R}|x-y|.$$

On the support ${\rm Supp} (\Gamma)$ of $\Gamma$, we can define the distance $d_\Gamma$ by
$$ \forall x,y\in {\rm Supp} (\Gamma),\quad d_\Gamma(x,y)=d_e(\Gamma(x),\Gamma(y)).$$
Hence there is an isometry, namely $\Gamma$, between the metric-measure space $({\rm Supp} (\Gamma),d_\Gamma,M )$ and the metric-measure space $({\rm Supp} (\chi),d_e, M(B_R)\times\lambda_R)$. Since the closure of ${\rm Supp} (\chi)$ with respect to the Euclidian distance is equal to $B_R$, the completion of the metric space $({\rm Supp} (\Gamma),d_\Gamma)$, denoted by $(C,d_\Gamma )$,   is isometric to $(B_R,d_e)$. That isometry, which coincides with $\Gamma$ on $ {\rm Supp} (\Gamma)$ is still denoted by $\Gamma$ and its inverse, which coincides with $\chi$ on $ {\rm Supp} (\chi)$, is still denoted by $\chi$. Obviously, the metric space $(C,d_\Gamma )$ is compact. 

Furthermore, since ${\rm Supp} (\Gamma)$ is a Borelian subset of $B_R$ as well as a Borelian subset of $C$ (for the respective topologies of $B_R$ and $C$) and since $M({\rm Supp} (\Gamma))=M(B_R)$, the measure $M$ can be extended to the whole of the Borelian subsets of $C$ by prescribing: 
$$\text{ for any Borelian subset A of C}, \quad M(A)=M(A\cap {\rm Supp} (\Gamma)).$$

Via pullback, the space $C$ inherits the structure of Riemannian manifold (smooth, complete, connected, $m$-dimensional manifold equipped with a smooth metric tensor). The (only) chart is given by $\Gamma:C\to B_R$. We summarize below what we have proved as well as the properties inherited from the pullback metric:

\begin{theorem}\label{thinflevy}
If $\psi'(1)<1$, we can find a compact Riemannian manifold $(C,g)$ and a Borelian subset $B$ of $B_R$ such that:
\begin{enumerate}
\item $B$ is dense in $B_R$ for the Euclidian distance and has full $M$-measure, namely  $M(B_R\setminus B)=0$,
\item $C$ is the completion of $B$ with respect to the geodesic distance on $C$,
\item the volume form on $C$ coincides with the measure $M$ on $B$,
\item in the system of local coordinates given by the chart $\Gamma$,  the Riemannian metric tensor on $C$ reads $$g=d(\omega)(dx_1^2+\dots+dx_m^2)\quad \text{ with }\quad d(\omega)=\frac{M(B_R)^2}{C_R^2}.$$
\end{enumerate}
\end{theorem}

%%%%%%%%%%%%%%%%%%%%%%%%%%%%%%%%%%%%%%%%%%%%%%%%%%%%%%%%%%%%%%%%%%%%%%%%
 \subsubsection{The strongly intermittent case $1 \leq \psi'(1)<m$}
%%%%%%%%%%%%%%%%%%%%%%%%%%%%%%%%%%%%%%%%%%%%%%%%%%%%%%%%%%%%%%%%%%%%%%%%

Now we consider a log-infinitely divisible random measure $M$  satisfying $1\leq \psi'(1)<m$. The measure $M$ gives mass to small sets so that we cannot use classical theorems of optimal transport theory if we consider the ambiant space $B_R$ equipped with its Euclidian structure. It thus seems hopeless to solve the problem:
\begin{align}\label{optiME}
 \text{Find the mapping }\varphi  \text{ realizing the infimum: }
\inf_{\substack{\varphi:B_R\to B_R \\\varphi_{\#\overline{M}}=\lambda_R}}\int_{B_R}|\varphi(x)-x|^2\overline{M}(dx).
\end{align}

Our strategy is the following: can we find a pair $(M',T')$, where $M'$ is a probability measure and $T'$ is an optimal (for the Euclidian quadratic cost) transport map pushing $M'$ forward to $\lambda_R$, such that the following Monge type optimization problem possesses a (unique) solution?
\begin{align}\label{optiM}
 \text{Find the mapping }\varphi  \text{ realizing the infimum: }
\inf_{\substack{\varphi:B_R\to B_R \\\varphi_{\#\overline{M}}=\lambda_R}}\int_{B_R}|\varphi(x)-T'(x)|^2\overline{M}(dx).
\end{align}
It turns out that the above problem is a mathematical formulation of the following intuitive argument: if we cannot find an optimal transport pushing $\overline{M}$ forward to $\lambda_R$, can we find an intermediate measure $M'$ and optimal transports $T,T'$ respectively pushing $\overline{M}$ forward to $M'$ and $M'$ forward to $\lambda_R$? If the answer is positive, then the composition $T'\circ T$ pushes $M$ forward to $\lambda_R$. Though the composition $T'\circ T$ is in general not optimal for the Euclidian quadratic cost, it is the composition of two gradients of convex functions, which is not  bad in terms of regularity. So we have formalized some kind of two-step optimal transport. And more generally, if we cannot find a two-step optimal transport, is it possible to find a $n$-step optimal transport between $\overline{M}$ and $\lambda_R$, that is a composition of $n$ gradients of convex functions?

The reader may have the following objection: from classical theorems, existence of a (unique) optimal transport map pushing $\overline{M}$ forward to $M'$ or $\lambda_R$ does not depend on the target measure ($M'$ or $\lambda_R$) but only on $\overline{M}$ through the fact that $\overline{M}$ does or does not give mass to small sets. So, a priori, two-step optimal transports may be as difficult to exhibit as optimal transports. Our idea lies in the fact that Problem \eqref{optiM} can be reformulated in a very simple way if we equip the ball $B_R$ with an appropriate Riemannian structure. Indeed, if we change the unknown in \eqref{optiM}, we get the following equivalent problem
\begin{align}\label{optiM2}
 \text{Find the mapping }\psi  \text{ realizing the infimum: }
\inf_{\substack{\psi:B_R\to B_R \\\psi_{\#\overline{M}}=M'}}\int_{B_R}|T'(\psi(x))-T'(x)|^2\overline{M}(dx).
\end{align}
It turns out that we can equip the ball $B_R$ with a Riemannian structure, the distance of which matches $d(x,y)=|T'(x)-T'(y)|$ and the volume form of which matches $M'$ (up to a multiplicative constant). Problem \eqref{optiM2} thus reduces to a classical problem of optimal transportation theory on smooth Riemaniannian manifolds: for such a mapping to exist, it is mainly sufficient that $\overline{M}$ does not give mass to the small sets associated to the distance $d$. That is the main constraint when choosing the measure $M'$: it must be the volume form associated to a Riemannian metric, the small sets of which are not charged by $M$. Of course, the argument generalizes to $n$-step optimal transports.

The main difficulty thus lies in choosing the number $n$ of steps  and the intermediate measures. The crucial point is the following. Given $n\geq 1$,  $M$ can be seen as the composition of $n$ multiplicative chaos (see subsection \ref{compo}): we can find $n$ independent independently scattered log-infinitely divisible random measures $\mu^{(1)}, \dots,\mu^{(n)}$  associated to $(\varphi/n,\theta)$ (see subsection  \ref{cons:loglevy}). The corresponding processes $\omega_l$ associated to $\mu^{(1)}, \dots,\mu^{(n)}$  are respectively denoted by $\omega_l^{(1)},\dots,\omega_l^{(n)}$.  We define recursively for $k\leq n$:
 $$M^{(0)}(dx)=dx \quad \text{and }\quad M^{(k)}(dx)=\lim_{l\to 0}\,\,e^{\omega_l^{(k)}(x)}M^{(k-1)}(dx),$$ where the limits have to be understood in the sense of weak convergence of Radon measures. Then both measures $M$ and $M^{(n)}$ have the same law so that we assume, in what follows, that $M$ and $M^{(n)}$ coincide.
%Let us choose $n$ such that
% \begin{equation}\label{defn}
% d\psi(2)<n(d-\psi'(1)),%+\psi'(1). 
%\end{equation}
%so as to make valid the relation
%\begin{equation}\label{defkn}
%\forall k\leq n-1,\quad \frac{d\psi(2)}{n}\leq d-\frac{k}{n}\psi'(1).
%\end{equation}
%From \eqref{defkn}, the Laplace exponent $\psi^{(1)}=$ of the measure $M^{(1)}$ satisfies  $$d-\psi^{(1)'}(1)=d-\frac{1}{n}\psi'(1)>\frac{\psi(2)}{n}.$$
%That condition is sufficient to ensure that $M^{(1)}$ does not give mass to the Euclidian small sets (see section \ref{cf2}).
That procedure allows to see each measure $M^{(k)}$ as a chaos with respect to $M^{(k-1)}$ with a reduced intermittency parameter $\psi'(1)/n$. If we can equip the ball $B_R$ with a Riemannian metric $g^{(k-1)}$ the volume form of which coincides with $M^{(k-1)}$, it turns out that $M^{(k)}$ does not give mass to the $g^{(k-1)}$-small sets provided that the intermittency parameter $\psi'(1)/n$ is small enough. So it suffices to choose $n$ small enough. In consequence there is a unique optimal transport map (w.r.t. to the quadratic cost function associated to the metric $g^{(k-1)}$)  that pushes the renormalized measure $\overline{M}^{(k)}$ forward to the renormalized measure $\overline{M}^{(k-1)}$. And so on for the different values of $k\leq n$. Thus we claim:

\begin{theorem}\label{thsuplevy}
If $1\leq \psi'(1)<d$, we can find $n\geq 1$ and $n$ gradients of convex functions $T^{(1)},\dots,T^{(n)} $ such that $\forall k=1,\dots,n$, the mapping $\varphi^{(k)}=T^{(1)}\circ \cdots\circ T^{(k)}$ pushes the measure $\overline{M}^{(k)}$ forward to the Lebesgue measure and minimizes the quantity
$$\inf_{\substack{T:B_R\to B_R\\T_{\#\overline{M}^{(k)}}=\lambda_R}}\int_{B_R}|T(x)-\varphi^{(k-1)}(x)|^2\overline{M}^{(k)}(dx).$$
\end{theorem}
As a corollary, we get

\begin{theorem}\label{thsuplevyR}
We can find a compact Riemannian manifold $(C,g)$ and a Borelian subset $B$ of $B_R$ such that:
\begin{enumerate}
\item $B$ is dense in $B_R$ for the Euclidian distance and has full $M$-measure, namely  $M(B_R\setminus B)=0$,
\item $C$ is the completion of $B$ with respect to the geodesic distance on $C$.
\item the volume form on $C$ coincides with the measure $M$ on $B$
\item in the system of local coordinates given by the chart $\varphi^{(n)}$,  the Riemannian metric tensor on $C$ reads $$g=d(\omega)(dx_1^2+\cdots +dx_m^2)\quad \text{ with }\quad d(\omega)=\frac{M(B_R)^2}{C_R^2}.$$
\end{enumerate}
\end{theorem}

\begin{remark}
In Theorems \ref{thinflevy} and \ref{thsuplevy}, if $R<T$, by scale invariance, the random variable $d$ has the same law as $e^{2\Omega_{\frac{R}{T}}}\frac{M(B_T)^2}{C_T^2}$ where $\Omega_{\frac{R}{T}}$ is an infinitely divisible random variables the law of which is characterized by $\E[e^{q\Omega_{\frac{R}{T}}}]=\big(\frac{R}{T}\big)^{-\psi(q)}$. So the radius $R$ of the ball $B_R$ influences the random metric through a log-infinitely divisible random variable, which turns out to be log-normal for log-normal MRM.
%
%We can generalize our approach to more general sets $A$ (with, say, $C^1$-boundary) instead of the ball $B_R$. The random variable $d$ in \eqref{cov} is a function of the set $A$. In the case where the Gaussian distribution appearing in \eqref{cov} is the Gaussian Free Field (see also subsection \ref{app:KPZ} below), the law of $d$ is invariant under conformal transformations of $A$ since the law of the GFF is.
\end{remark}

%%%%%%%%%%%%%%%%%%%%%%%%%%%%%%%%%%%%%%%%%%%%%%%%%%%%%%%%%%%%%%%%%%%%%%%%

\FloatBarrier

\section{Applications}\label{sec:appl}
%%%%%%%%%%%%%%%%%%%%%%%%%%%%%%%%%%%%%%%%%%%%%%%%%%%%%%%%%%%%%%%%%%%%%%%%%%%%%%%%%%%%%%%%%%

%If we denote by $\pi_n$ the probability measure $({\rm Id},\chi_n)_\#\lambda_R$. This probability measure on $B_R\times B_R$ weakly converges as $n\to +\infty$ towards an optimal transference plan $\pi $ between  $\lambda_R$ and $\overline{M}$. Roughly speaking, one may convince oneself that the curve $(t,\gamma_n^{x,y}(t))_{0\leq t \leq 1}$ should converge towards a possibly multivariate curve lying in the support of $\pi$. In case $\pi$ is deterministic, that support is the graph of the gradient of some convex function $\Gamma$ so that the curve $(t,\gamma_n^{x,y}(t))_{0\leq t \leq 1}$ should converge towards the curve $$\Big(t,\chi\big((1-t)\Gamma(x)+t\Gamma(y)\big)\Big)_{0\leq t \leq 1}.$$ Hence there is much to learn of the convergence of geodesics from the existence of a unique transport map between $\overline{M}$ and $\lambda_R$. This is the purpose of our paper. Because of the singularity of the measures $\overline{M}$ and $\lambda_R$, the function $\Gamma$ has very bad regularity properties. Hence our purpose it to give a sufficient condition on the parameter $ \gamma$ of the measure $\overline{M}$ in order to have a unique optimal transport map between $\overline{M}$ and $\lambda_R$ and to give a meaning to the convergence of the curves.

\subsection{Application to the KPZ framework}\label{app:KPZ}
%%%%%%%%%%%%%%%%%%%%%%%%%%%%%%

A special case of \eqref{measureM} in dimension $2$ has recently received much attention. When $X$ is a Gaussian Free Field (GFF), that is a Gaussian process with covariance function given by the Green function of the Laplacian, the measure $M$ is called quantum measure. That situation is motivated by  Liouville quantum gravity \cite{cf:DuSh}.  Much effort has recently been made to understand the geometry of these measures. The first important step was to prove the so-called KPZ formula \cite{Benj,cf:DuSh,cf:RhoVar}. Roughly speaking, the KPZ formula gives the correspondence between the Haudorff dimension of a set as
seen by the Lebesgue measure and the the Haudorff dimension of this set as
seen by the measure $M$. More precisely, for a given compact set $E\subset B_R$ and a mesure $\nu$, the $s$-dimensional Hausdorff measure of $E$ wrt to $\nu$ is the quantity:
\begin{align}\label{intro:hausdorff}
H^s(\nu,E)&=\lim_{\delta\to 0}H^s_\delta(E,\nu),\quad  \text{ where }\\ H^s_\delta(\nu,E)&=\inf \{\sum_n \nu(B_n)^{s/2};E\subset \bigcup B_n, B_n\text{ measurable set with }0<M(B_n)\leq \delta\}\nonumber.
\end{align}
Then we define the Hausdorff dimension ${\rm dim}^\nu_H(E) $ of the set $E$ wrt to the measure $ \nu$: 
\begin{equation}\label{intro:dim}
{\rm dim}^\nu_H(E)=\inf\{s> 0;H^s(\nu,E)=0\}=\sup\{s>0;H^s(\nu,E)=+\infty\}.
\end{equation}
When $\nu$ is given by the Lebesgue measure, the corresponding Hausdorff measures and dimensions will be called Euclidian. The KPZ formula asserts that the Euclidian Hausdorff dimension ${\rm dim}_H(E)$  and the random Hausdorff dimension ${\rm dim}^M_H(E) $ are linked by the relation 
\begin{KPZ}
Almost surely, we have 
\begin{equation}\label{eqKPZ}
\xi\Big(\frac{{\rm dim}^M_H(E)}{2}\Big)={\rm dim}_H(E) 
\end{equation}
 where $\xi$ is the so-called structure exponent of $M$:
$$\xi(q)=\big(2+\frac{\gamma^2}{2}\big)q-\frac{\gamma^2}{2}q^2. $$ 
\end{KPZ}
Actually, Formula \eqref{eqKPZ} remains valid for any MRM as constructed in \eqref{cons:loglevy} as soon as $M$ possesses moments of negative order (see \cite{cf:RhoVar}).

A further understanding of the geometry of the measure $M$ is to define a Riemannian metric  such that the corresponding volume form coincides with $M$. Such a structure permits to define fundamental objects such as distance, arclength, geodesics that are naturally associated to the measure $M$. This is an open problem 
suggested by K. Izyurov at the conference "Les Diablerets". 

Theorem \ref{thsuplevyR} allows to understand the measure $M$ on $B_R$  as the volume form of a  Riemannian manifold (up to a set of null $M$-measure). The geodesics are easily described via the  transport maps.  By sticking to the notations of Theorem \ref{thsuplevy}, we define $\Gamma=\varphi^{(n)}$ and $\chi=\Gamma^{-1}$. Since $\Gamma$ is an isometry of metric spaces, a curve $\gamma:[0,1]\to C$ is a geodesic on $C$ if and only if $\Gamma(\gamma)\subset B_R$ is a geodesic for the Euclidian metric, that is a segment. The geodesic joining $x,y\in B\subset B_R$ (parameterized by constant speed)   is thus given by 
$$\gamma^{x,y}:t\in[0,1]\to \chi\big(t\Gamma(x)+(1-t)\Gamma(y)\big)\in C.$$
\begin{figure}
\centering 
%%----start of first subfigure---- 
\subfloat[Conformal geodesic $\gamma^2\simeq 2$]{
\label{geod:a} %% label for first subfigure 
\includegraphics[width=12.3cm,height=10cm]{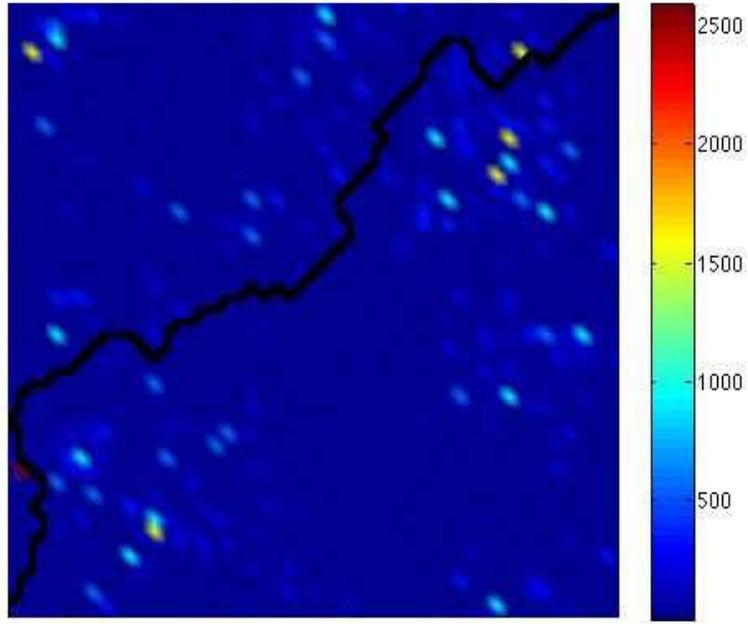}} 
\hspace{0.1\linewidth} 
%----start of second subfigure---- 
\subfloat[Comparison of conformal/optimal transport geodesics for $\gamma^2\simeq 1.5$]{ 
\label{geod:b} %% label for second subfigure 
\includegraphics[width=12.3cm,height=10cm]{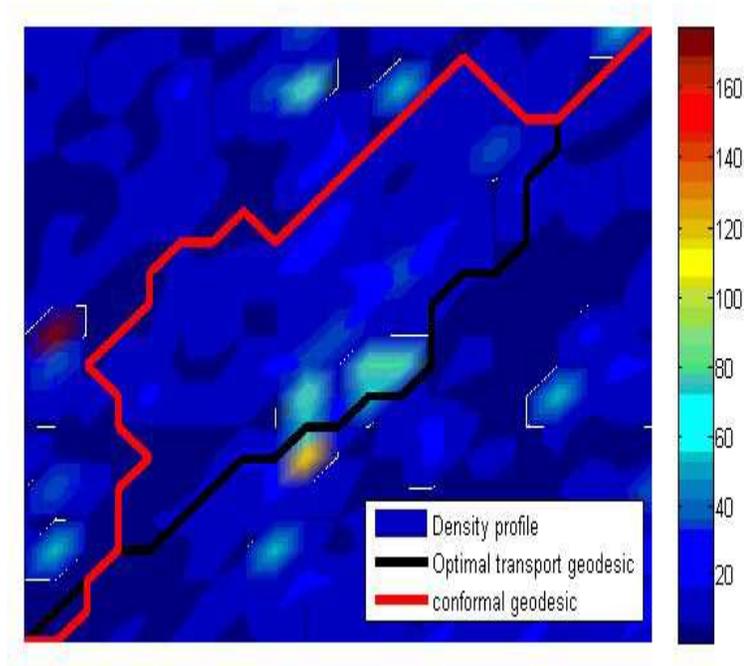}}
\caption{Simulations of conformal geodesics (weighted by the density of the MRM along the path) and geodesics obtained via optimal transport for a log-normal MRM.}
\label{geod} %% label for entire figure 
\end{figure}

In what follows, we give a formal interpretation of the KPZ formula as well as a conjecture on the Hausdorff dimension of the geodesics. Consider a compact set $K$, contained in ${\rm Supp}(\chi)$. Then $\chi(K)$ is a subset of $C$ but may also be thought of as a subset of $B_R$. Thus we can define the Hausdorff dimension of $\chi(K)$ as a subset of $C$ equipped with the distance $d$, denoted by ${\rm dim}_H^d(\chi(K))$ and   the Hausdorff dimension of $\chi(K)$ as a subset of $B_R$ equipped with the Euclidian distance, denoted by ${\rm dim}_H^e(\chi(K))$. By isometry between $(C,d,M)$ and $(B_R,d_e,M(B_R)\lambda_R)$, we  have the obvious relation $${\rm dim}_H^d(\chi(K))={\rm dim}_H(K).$$ Furthermore,  if we assume that formula \eqref{eqKPZ} still holds for these Hausdorff dimensions (defined in terms of distances not measures), we will have
$${\rm dim}_H^e(\chi(K))=\xi\Big(\frac{{\rm dim}^d_H(\chi(K))}{2}\Big)=\xi\Big(\frac{{\rm dim}_H^e(K)}{2}\Big).$$
The KPZ formula can thus be reinterpreted in terms of the graph of the function $\chi$: $\chi$ sends a set $K$ with Hausdorff dimension $q$ towards the set $\chi(K) $ with Hausdorff dimension $\xi(q/2)$. In particular, the geodesics with respect to the metric $d$ should have Hausdorff dimension $\xi(1/2)=1+\frac{\gamma^2}{8}$.

We cannot give a rigorous argument to that conjecture because  the general theory of optimal transport suffers a definite lack of  strong estimates concerning the optimal transport maps. So we cannot prove the KPZ formula connecting the Hausdorff dimensions defined in terms of distances (not measures) and we leave that point as an open question.

%%%%%%%%%%%%%%%%%%%%%%%%%%%%%%%%%%%%%%%%%%%%%%%%%
\subsection{Multifractal random change of time}

Another important motivation for constructing an optimal transportation maps associated with MRMs is the possibility of defining multifractal random change of times as suggested by Mandelbrot in the 1d case: if $B$ is a 1-d Brownian motion and $M$ is a 1-d MRM, we can define the process $t\mapsto B_{M([0,t])} $. This process has been widely used in financial modelling since it enjoys nice properties: it is a square-integrable martingale with long-range correlated stationary increments and power law spectrum. Extending that construction may be done via optimal transportation maps.

Let us fix $R>0$ and a $m$-dimensional MRM as constructed in \ref{cons:loglevy}. Let $\Gamma,\chi$ be the 
maps pushing $M$ forward to $\frac{M(B_R)}{C_R}\, dx$ and vice versa (sticking to the previous notations, we denote by $B$ the support of $\Gamma$). We further consider a Gaussian white noise $W(dx)$ defined on $\R^m$. Given $x\in \R^m$, we denote by $C(x)$ the cube $[0,x_1]\times\cdots \times [0,x_m]$ (with the convention that, if $x_i<0$, the interval $[0,x_i]$ stands for $[x_i,0]$). Finally, we define the process 
$$\forall x\in B_R,\quad B(x)=W\big(\Gamma(B\cap C(x))\big).$$
It is readily seen that the process $B$ is isotrop and stochastically scale invariant:
$$\forall \lambda \leq 1,\quad \big(B( \lambda x)\big)_{x\in B_R}\stackrel{law}{=}\lambda^{m/2}e^{\frac{1}{2}\Omega_\lambda} \big(B( x)\big)_{x\in B_R},$$ where $\Omega_\lambda$ is an infinitely divisible random variable, independent of $ \big(B( x)\big)_{x\in B_R}$, the law of which is characterized by:
$$ \E[e^{iq \Omega_\lambda}]=\lambda^{-\varphi(q)}.$$

%Furthermore, the Kolmogorov criterion ensures the almost sure continuity of the process.
Of course, we can generalize this approach to make a multifractal random changes of time to many other multidimensional stochastic processes like Levy noise, fractional Brownian motion, and so on... 
\section{Proofs}
%%%%%%%%%%%%%%%%%%%%%%%%%%%%%%%%%%%%%%%%%%%%%%%%%%%%%%%%%%%
We first prove the key estimate of the paper.
%%%%%%%%%%%%%%%%%%%%%%%%%%%%%%
\subsection{Fundamental result}\label{fund}
%%%%%%%%%%%%%%%%%%%%%%%%%

Let $B$ be a Borelian subset of $B_R$ and $\kappa$ be a probability measure on $B_R$ supported by $B$, meaning $\kappa(B)=1$. We assume that the set $B$ is equipped with a distance $d$ and that  the completion of $B$ with respect  to the distance $d$, denoted by $(C,d)$, is compact. We assume  that the Borelian subsets of $B$ with respect to the Euclidian topology coincide with the Borelian subsets of $B$ w.r.t. the distance $d$  so that we can extend the measure $\kappa$ to the whole Borelian sigma-algebra of $C$ by prescribing
$$\forall A \subset C \,\text{ Borelian}, \quad \kappa(A)=\kappa(A\cap B).$$

{\bf The classes $\mathbf{R_\alpha,R_\alpha^-}$.} For any $\alpha>0$, we introduce the set $R_\alpha$ of Radon measures $\nu$ on $C$ satisfying:   for any $\epsilon>0$, there are $\delta>0$, $D>0$ and a $d$-compact subset  $ K_\epsilon\subset C$ such that $\nu(C\setminus K_\epsilon)<\epsilon $ and the measure $\nu_{\epsilon}(dx)=\one_{K_\epsilon}(x)\nu(dx)$ satisfies
\begin{equation}\label{holderd}
\forall U\,d\text{-open ball},\quad \nu_{ \epsilon}(U)\leq D\times\,{\rm Diam}_d(U)^{\alpha+\delta},
\end{equation}
where ${\rm Diam}_d$ denotes the diameter of the ball $U$ with respect to the distance $d$. We further define the set of Radon measures $R^-_\alpha=\bigcap_{\beta<\alpha}R_\beta$.
%We denote by $B_d(x,r)$ the ball of center $x$ and radius $r$ with respect to $d$. 
%The measure $M'$ is obviously in the class $R'_2$ since
%\begin{align*}
%M'(B_\Gamma(x,r))=&\int\one_{B(\Gamma(x),r)}(\Gamma(y))\overline{M}'(dy)\\
%=&\int\one_{B(\Gamma(x),r)}(y)\lambda_R(dy)\\
%=&Dr^2
%\end{align*}
%where $D$ is the Lebesgue measure of the unit Euclidian ball.

We further give an energy condition for  a Radon measure $\nu$ to be in the class $R_\alpha^-$.
For $\alpha>0$, we define
\begin{equation}\label{energie}
C_\alpha(\nu)=\int_{C}\int_{C}\frac{1}{d(x,y)^\alpha}\nu(dx)\nu(dy).
\end{equation}
It is plain to see that 
$$\nu\text{ satisfies }C_\alpha(\nu)<+\infty\Rightarrow \nu\in R_\alpha^-.$$
Conversely a measure $\nu$ obeying \eqref{holderd} satisfies $C_\beta(\nu)<+\infty$ for each $\beta<\alpha+\delta$.

{\it Notations:} In what follows, we use the superscript $e$ (i.e. $R_\alpha^e$, $R_\alpha^{e-} $ and $C_\beta^e$) to mention that the  distance $d$ is equal to the Euclidian distance.\qed

\begin{proposition}\label{fundprop}
Assume that the measure $\kappa$ belongs to $R^e_\alpha$.
Let $N$ be the Radon measure on $B_R$ defined, $\P$-almost surely, as the following limit (in the sense of weak limit of Radon measures):
$$N(dx)=\lim_{l\to 0}N_l(dx),\quad \text{ where }N_l(dx)\stackrel{\rm{def.}}{=}\,\,e^{\omega_l(x)}\kappa(dx).$$
If $\psi(2)<\alpha$ then  the martingale $\big(N_l(B)\big)_l$ is uniformly integrable and, consequently, $N$ is non-trivial and almost surely supported by $B$.
\end{proposition}

\noindent {\it Proof.} We remind the reader of the fact that the family $\big(N_l(B)\big)_l$ is a right-continuous positive martingale and thus converges almost surely.

 Since $\kappa$ belongs to $R_\alpha$, for each $\epsilon>0$, we can find $\delta>0$, $D>0$ and a compact subset $K$ such that $\kappa(C\setminus K)<\epsilon$ and 
$$\limsup_{r\to \infty}\frac{\ln \kappa(B_r^x)}{r}<-\alpha-\delta\quad \text{uniformly for }x\in K,$$ where $B_r^x$ stands for the $d$-ball of radius $e^{-r}$ and center $x$.
This implies $C_\beta(\kappa_K)<+\infty$ for each $\alpha<\beta<\alpha+\delta$.

Then we compute
$$ \E\Big[N_l(B\cap K)^2\Big]=\int_{B\cap K}\int_{B\cap K}e^{\psi(2)K_l(|y-x|)}\kappa(dx)\kappa(dy)$$ where $K_l$ is given on $\R^m$ by $$K_l(x)=\int_{G}\rho_l(gx)\,H(dg)$$
where $\rho_l(y)$ ($y\in\R$) is defined by 
$\ln(T/|y|)$ if $l\leq |y|\leq T$, $\rho_l(y)=\ln(T/l)+1-|y|/l$ if $|x|\leq l$ and $0$ otherwise. In particular, on $B$, $K_l$ is not greater than $C+\ln\frac{T}{|x|} $ for some positive constant $C$. As a consequence, for some positive constant $C'$ which may change along the inequalities, we have
$$\sup_l \E\Big[N_l(B\cap K)^2\Big]\leq C'\int_{B\cap K}\int_{B\cap K}e^{\psi(2)\ln\frac{T}{|y-x|}}\kappa(dx)\kappa(dy)\leq C' C_\beta(\kappa)<+\infty$$ since $\psi(2)< \beta$. The martingale $\big(\widetilde{\nu}_l(B\cap K)\big)_l$ is bounded in $L^2(\Omega)$ and is therefore uniformly integrable. It is plain to deduce that the martingale $\big(\widetilde{\nu}_l(B)\big)_l$ is  uniformly integrable.\qed

\begin{theorem}\label{fundth}
1) Assume that the measure $\kappa$ belongs to $R_\alpha\cap R^e_\varsigma$  for some $\varsigma>\psi(2)$. Then    $$N\in R_{\alpha\frac{\varsigma-\psi(2)}{\varsigma}+\psi(2)-\psi'(1)}.$$ Consequently, we also have: $$\kappa\in R^{-}_\alpha\cap R^{e-}_\varsigma\text{ for some }\varsigma>\psi(2)\,\,\,\Rightarrow\,\,\, N\in R^-_{\alpha\frac{\varsigma-\psi(2)}{\varsigma}+\psi(2)-\psi'(1)}.$$

2) In particular, we have in the Euclidian case: if $\kappa\in R_\alpha^e$ (resp. $R_\alpha^{e-} $) and $\psi(2)<\alpha$ then $M\in R^{e}_{\alpha-\psi'(1)}$ (resp. $R^{e-}_{\alpha-\psi'(1)}$).
\end{theorem}

\noindent {\it Proof.}  %By convexity arguments, note that $\psi(2)<\varsigma$ implies $\psi'(1)<\varsigma$. 
Since $\kappa$ belongs to $R_\alpha\cap R_\varsigma^e$, for each $\epsilon>0$, we can find $\delta>0$, $D>0$ and a compact subset $K'$ for the $d$-topology ($d$-compact for short) such that $\kappa(C\setminus K')<\epsilon/2$ and 
\begin{equation}\label{limsup}
\limsup_{r\to \infty}\frac{\ln \kappa(B_r^x)}{r}<-\alpha-\delta\quad \text{uniformly for }x\in K',
\end{equation}
 where $B_r^x$ still stands for the $d$-ball of radius $e^{-r}$.

Since $\kappa$ belongs to $R^e_\varsigma$, for each $\epsilon>0$, we can find  a compact subset $K''$ for the Euclidian topology such that $\kappa(B_R\setminus K'')<\epsilon/4$ and  $C_\beta(\one_{K''}(x)\kappa(dx))<+\infty$ for each $\beta<\varsigma$. $K''\cap B$ is a Borelian subset of $B$ (for the Euclidian topology) and hence a Borelian subset of $C$. Since $M$ is inner regular on $C$ (recall that $C$ is compact), we can find a $d$-compact set $\widetilde{K}$ of $C$ such that $\widetilde{K}\subset K''$ and $M(K''\setminus \widetilde{K})<\epsilon/4$. Finally we set $K=K'\cap \widetilde{K}$, which is $d$-compact. The set $K$ also satisfies: $M(C\setminus K)<\epsilon$, \eqref{limsup} is valid on $K$ and $C^e_\beta(\one_{B\cap K})<+\infty
$ for any $\beta\leq \varsigma$.

Even if it means multiplying $\kappa$ by a constant, we assume $\kappa(K)=1$. We  consider on $\Omega\times K$ the probability measure $\Q$ defined by
$$\int_{\Omega\times K}f(\omega,x)\,d\Q=\E\Big[\int_{K}f(\omega,x)N(dx)\Big]$$ for all measurable nonnegative functions $f$.
  
Given $l'<l$, we define $$\omega_{l',l}(x)=\omega_l(x)-\omega_{l'}(x).$$ For a sequence $l_1<\cdots<l_n$,  the random variables $\omega_{l_1,l_2},\dots,\omega_{l_{n-1},l_n}$ are $\Q$-independent and 
$$\int e^{\lambda\omega_{l_i,l_j}}\,d\Q=\E\big[e^{(1+\lambda)\omega_{l_i,l_j}(x)}\big]=e^{\psi(1+\lambda)\ln(l_j/l_i)}.$$
The process $u\in\R_+\mapsto \omega_{e^{-u}}$ is therefore an integrable L\'evy process (we can consider a version that is right-continuous with left limits). From the strong law of large numbers, we have
$$\Q\text{ a.s.},\quad \frac{\omega_{e^{-u}}}{u}\to \psi'(1)\quad \text{as }u\to \infty .$$% where $a_n=\sum_{p=1}^nk_p(0)=\gamma^2\ln(nT)+2\gamma^2$.
This implies that, $\P$ a.s.,
$$N \text{ a.s.}, \quad \frac{\omega_{e^{-u}}}{u}\to \psi'(1)\quad \text{as }u\to \infty.  $$
Therefore, $\P$ a.s., for each $\epsilon>0$ we can find a compact $K^1_\epsilon\subset K$ such that $N(K\setminus K^1_\epsilon)<\epsilon$ and  $\frac{\omega_{e^{-u}}(x)}{u}\to \psi'(1)$ uniformly w.r.t. $x\in K^1_\epsilon$ as $u\to \infty$. Now we define 
$$N_q(dy)=\lim_{l\to 0}e^{\omega_{l,e^{-q}}(y)}\,\kappa(dy)\quad \text{and}\quad P_q(x)=\int_{B^q_x\cap K}N_q(dy).$$ We further define the function $\theta_q$  by 
$$\theta_q(x,y)=\left\{\begin{array}{ll}1 &\text{ if }d(x,y)\leq e^{-q}\\
  0 &\text{otherwise.}\end{array}\right.$$ Thus we have
$P_q(x)=\int_{K}\theta_{q}(x,y)N_q(dy)$ and
\begin{align*}
\int P_q\,d\Q&=\E\int_{K}\int_{K}\theta_q(x,y)N_q(dy)N(dx)\\
&\leq \int_{B\cap K}\int_{B\cap K}\theta_q(x,y)e^{\psi(2)\big(C+\ln(\frac{e^{-q}}{|x-y|})\big)}\,\kappa(dx)\kappa(dy)\\
%&=\int_{B_R}\int_{B_R}\chi_q(x,y)e^{-r\psi(2)}\frac{1}{|t-s|^{\psi(2)}}\nu(ds)\nu(dt),
%&\leq \int_{B_R}\int_{B_R}\theta_q(x,y)e^{\phi(2)\ln\frac{1}{q|y-x|}}\,M'(dx)M'(dy)\\
&\leq e^{\psi(2)C}\int_{B\cap K}\int_{B\cap K}\theta_q(x,y)e^{-q\psi(2)}  \frac{1}{|y-x|^{\psi(2)}}\,\kappa(dx)\kappa(dy).
\end{align*}
By using the above relation, we obtain
\begin{align*}
\int\sum_{n\geq 1}e^{\beta  n}P_{n}\,d\Q& =\sum_{n\geq 1}\int_{B\cap K}\int_{B\cap K}e^{(\beta-\psi(2)) n}\chi_{n}(y,x)\frac{1}{|y-x|^{\psi(2)}}\,\kappa(dx)\kappa(dy).
\end{align*}
Note that (for some positive constant $D$)
$$\sum_{n\geq 1}e^{(\beta-\psi(2)) n}\chi_{n}(y,x) =\sum_{1\leq n\leq-\ln d(x,y)}e^{(\beta-\psi(2)) n} \leq D\frac{1}{d(x,y)^{\beta-\psi(2)}}$$ in such a way that we obtain
\begin{align*} 
 \int\sum_{n\geq 1}e^{\beta  n}P_{n}\,d\Q&\leq D\int_{B\cap K}\int_{B\cap K}\frac{1}{|y-x|^{\psi(2)}}\frac{1}{d(x,y)^{\beta-\psi(2)}}\,\kappa(dx)\kappa(dy)
 %\\&= D\int_{B_R}\int_{B_R}\frac{1}{|\chi(y)-\chi(x)|^{\phi(2)}}\frac{1}{|y-x|^{\beta-\phi(2)}}\,dxdy.
\end{align*} 
We want to prove that the latter integral is finite for a well chosen $\beta$.  We fix $$\psi(2)<\beta=(\alpha+\delta/2)\frac{\varsigma-\psi(2)}{\varsigma}+\psi(2).$$
 We consider  $p,q>1$ satisfying the relation $\frac{1}{p}+\frac{1}{q}=1$ and given by 
  $$p=\frac{\varsigma}{\psi(2)}\quad\text{ and }\quad q=\frac{\varsigma}{\varsigma-\psi(2)}.$$ By H\"older's inequality, we have
\begin{align*}
\Big[\int_{B\cap K}&\int_{B\cap K}\frac{1}{|y-x|^{\psi(2)}}\frac{1}{d(x,y)^{\beta-\psi(2)}}\,\kappa(dx)\kappa(dy)\Big]\\
&\leq \Big[\int_{B\cap K}\int_{B\cap K}\frac{1}{|y-x|^{p\psi(2)}}\,\kappa(dx)\kappa(dy)\Big]^{1/p}\Big[\int_{B\cap K}\int_{B\cap K}\frac{1}{d(x,y)^{q(\beta-\psi(2))}}\,\kappa(dx)\kappa(dy)\Big]^{1/q}\\
&\leq \Big(C^e_{\varsigma}(\one_{K\cap B}(x)\kappa(dx))\Big)^{1/p}\Big(C_{\alpha+\delta/2}(\one_{K\cap B}(x)\kappa(dx))\Big)^{1/q},
\end{align*}
in such a way that the above integrals are finite.

 We deduce that, $\Q$ a.s., $e^{\beta n}P_{n}\to 0$ as $n\to \infty$. Therefore, $\P$ a.s., $e^{\beta n}P_{n}\to 0$ as $n\to \infty$  $ N$-almost surely. So we can find a compact $K^2_\epsilon\subset C$ such that $N(C\setminus K^2_\epsilon)<\epsilon$ and
 $$\limsup_{n\to \infty} \frac{\ln P_{n}(x)}{n}\leq -\beta\quad \text{uniformly for }x\in K^2_\epsilon.$$
 Finally we can set $\bar{K}=K^1_\epsilon\cap K^2_\epsilon$ and $N_{\bar{K}}(dx)= \one_{\bar{K}}(x)N(dx)$. We obtain
 \begin{align*}
 \limsup_{n\to \infty}\frac{ \ln N_{\bar{K}}(B_n^x)}{n}& =  \limsup_{n\to \infty}\frac{ \ln \int_{\bar{K}\cap B_n^x}e^{\omega_{e^{-n}}(x)}N_{n}(dx)}{n}\\
 &\leq -\beta+\psi'(1)
 \end{align*}
 uniformly w.r.t. $t\in \bar{K}$. We have proved $N\in R_{\alpha\frac{\varsigma-\psi(2)}{\varsigma}+\psi(2)-\psi'(1)}$ $\P$-almost surely. \qed
 
\subsection{Composition of MRM}\label{compo}
 %%%%%%%%%%%%%%%%%%%%%%%%%%%%%%%%%%%%%%%%%%%%%
 
 Now we prove that a non-degenerate log-infinitely divisible random measure $M$, i.e. satisfying $\psi(2)<+\infty$ and $\psi'(1)<d$, can be decomposed as an iterated multiplicative chaos. 
 
 Since $\psi(2)<+\infty$ we can find an integer $n$ such that
 \begin{equation}\label{defn}
 m\psi(2)<n(m-\psi'(1)).
\end{equation}
  Then we can find $n$ independent independently scattered log-infinitely divisible random measures $\mu^{(1)}, \dots,\mu^{(n)}$  associated to $(\varphi/n,\theta)$ (remind of the definition in subsection  \ref{cons:loglevy}). We assume that the random measures $\mu^{(1)}, \dots,\mu^{(n)}$ are constructed on the probability spaces $(\Omega^{(1)},\P^{(1)}),\dots,(\Omega^{(n)},\P^{(n)}) $. We define $\Omega=\Omega^{(1)}\times\dots\times\Omega^{(n)}$ equipped with the product $\sigma$-algebra and the product probability measure $\P=\P^{(1)} \otimes\dots\otimes \P^{(n)}$.   The corresponding processes $\omega_l$ associated to $\mu^{(1)}, \dots,\mu^{(n)}$  are respectively denoted by $\omega_l^{(1)},\dots,\omega_l^{(n)}$. Finally, we denote by $\E^{(i)}$ the conditional expectation given the variables $(\mu^{(k)})_{k\not=i}$.

 We define recursively for $k\leq n$:
 $$M^{(0)}(dx)=dx \quad \text{and }\quad M^{(k)}(dx)=\lim_{l\to 0}\,\,e^{\omega_l^{(k)}(x)}M^{(k-1)}(dx),$$ where the limits have to be understood in the sense of weak convergence of Radon measures.
Note that the choice of $n$ makes valid the relation
$$\forall k\leq n-1,\quad \frac{m\psi(2)}{n}<m-\frac{k}{n}\psi'(1).$$
Hence, we can apply recursively Proposition \ref{fundprop} and Theorem \ref{fundth} (with the distance $d$ equal to the Euclidian distance and $\kappa\in R^{e-}_2$ is the Lebesgue measure) to prove that, for each $k\leq n$, $$M^{(k)}\in R^{e-}_{2-\frac{k}{n}\psi'(1)}\quad \text{ and }\quad \E^{(k)}[M^{(k)}(B_R)]=M^{(k-1)}(B_R)\,\,\quad \P \text{ a.s.}.$$
 Thus we have $M^{(n)}\in R^{e-}_{m-\psi'(1)}$ and $\E[M^{(n)}(B_R)]=\lambda(B_R)$ (the Lebesgue measure of $B_R$).
 
What we now want to prove is  that the measure $M^{(n)}$ has the same law as the measure 
$$M(dx)=\lim_{l\to 0}\,\,e^{\omega_l^{(1)}(x)+\cdots+\omega_l^{(n)}(x)}\,dx.$$ 
We consider on $\Omega$ the $\sigma$-algebra $ \mathcal{G}_l$ generated by $\{\omega_r^{(1)}(x),\dots,\omega_r^{(n)}(x);x\in\R^m,T>r>l\}$. The conditional expectation of $M^{(n)}(E)$ w.r.t. $ \mathcal{G}_l$ is easily computed since, for each $k\leq n$, the martingale $(M_l^{(k)}(A))_l$ is $\P^{(k)}$-uniformly integrable. Indeed, we have:
 \begin{align*}
 \E\big[M^{(n)}(A)|\mathcal{G}_l\big]&=  \E\big[\E\big[M^{(n)}(A)|\mu^{(1)},\dots,\mu^{(n-1)},(\omega_r^{(n)}(x))_{x\in\R^m,T>r>l}\big]|\mathcal{G}_l\big]\\
 &=  \E\big[\E^{(n)}\big[M^{(n)}(A)|(\omega_r^{(n)}(x))_{x\in\R^m,T>r>l}\big]|\mathcal{G}_l\big]\\
 &=\E\big[\int_Ae^{\omega_l^{(n)}(x)}M^{(n-1)}(dx)|\mathcal{G}_l\big]\\
 &=\cdots\\
 &=\int_Ae^{\omega_l^{(n)}(x)+\cdots+\omega_l^{(1)}(x)} \,dx.
 \end{align*}
This latter quantity has the same law as $M_l(A)$. Since the martingale $\big(\E\big[M^{(n)}(A)|\mathcal{G}_l\big]\big)_l$ is uniformly integrable, we deduce that the family $(M_l(A))_l$ is uniformly integrable. Hence, both random variables $M(A)$ and $M^{(n)}(A)$ have the same law. In particular, $M\in R^e_{m-\psi'(1)}$. 
 
 \begin{corollary}\label{corocompo}
 If $\psi(2)<+\infty$ and $\psi'(1)<m$, then $M$ belongs to $R^{e-}_{m-\psi'(1)}$.
 \end{corollary}
 
 \begin{remark}
 The same composition argument shows that if a measure $\kappa\in R^e_\alpha$ and if $M$ is defined as the limit
 $$M(dx)=\lim_{l\to 0}\,\,e^{\omega_l(x)}\kappa(dx)$$ and satisfies $\psi(2)<+\infty$ then $M\in R_{\alpha-\psi'(1)}$. Though we won't use it in this paper, this is an important result concerning the structure of the support of  multifractal random measures.
   \end{remark}

 %%%%%%%%%%%%%%%%%%%%%%%%%%%%%%%%%%%%%%%%% 
\subsection{Proof of Theorem \ref{gamsmallevy}}\label{cf1}
 
 In the case where $\psi(2)<+\infty$ and $\psi'(1)<1$, we will prove that the measure $M$ does not give mass to small sets.  Let $S\subset B_R$ be a small set, that is a set with Hausdorff dimension (w.r.t. the Euclidian distance) not larger than $m-1$. 

From Corollary \ref{corocompo}, $M$ belongs to the class $R^{e-}_{m-\psi'(1)}$.   Since $\psi'(1)<1$, $M$ then belongs to the class $R^e_{m-1+\beta}$ for some $\beta>0$. We fix $\epsilon>0$. So, $\P$ a.s., we can find a compact set $K\subset B_R$ and $\delta,D>0$ such that $M(B_R\setminus K)\leq \epsilon$ and for all open balls $U\subset B_R$:
 $$M(U\cap K)\leq C{\rm diam}_e(U)^{m-1+\beta+\delta}.$$ 
 Since $m-1+\beta+\delta>{\rm dim}_H(S)$, we can find a covering of $S$ by open balls $(U_i)_i$ such that 
 $$\sum_i{\rm diam}_e(U_i)^{m-1+\beta+\delta}<\epsilon.$$
 Then we have
 $$M(S)\leq M(S\setminus K)+M(S\cap K)\leq \epsilon+\sum_iM(U_i\cap K)\leq \epsilon+\sum_i{\rm diam}_e(U_i)^{1+\beta+\gamma}\leq 2\epsilon.$$
 As we can make $\epsilon$ as small as we please, we deduce $M(S)=0$.
 \qed
 %%%%%%%%%%%%%%%%%%%%%%%%%%%%%%%%%%%%%%%%%%%%%%%
 
\subsection{Proof of Theorems \ref{thsuplevy} and \ref{thsuplevyR}}\label{cf2}
 %%%%%%%%%%%%%%%%%%%%%%%%%%%%%%%%%%%%%%%%%%%%%%%%%
 We proceed recursively to prove the existence of an optimal transport between the measure $M^{(k-1)}$ and $M^{(k)}$ on some appropriate Riemannian manifold:
 
 1) Step 1: We focus on the measure $M^{(1)}$, which has structure exponent $ \xi^{(1)}(q)=m q-\frac{\psi(q)}{n}$. Relation \eqref{defn} and the convexity of $\psi$ imply the following inequalities 
 $$\frac{\psi'(1)}{n}\leq \frac{\psi(2)}{n}=\frac{1}{m}\frac{m\psi(2)}{n}=\frac{1}{m}(m-\psi'(1))<1.$$ 
 Hence we can apply Theorem \ref{gamsmallevy} to find two optimal transport maps $\chi^{(1)},\Gamma^{(1)}$ that respectively push $\lambda_R$ forward to $\overline{M}^{(1)}$ and vice versa. Furthermore,  the quantity 
 $$\inf_{\substack{T:B_R\to B_R\\T_{\#}\overline{M}^{(1)}=\lambda_R}} \int_{B_R}|T(x)-x|^2\overline{M}^{(1)}(dx)$$ is achieved  at $\Gamma^{(1)}$.  We can also apply Theorem \ref{thinflevy} to find  a compact Riemannian manifold $(C^{(1)},g^{(1)})$ and a Borelian subset $B^{(1)}$ of $B_R$ such that:

-$B^{(1)}$ is dense in $B_R$ for the Euclidian distance and has full $M$-measure, that is  $M^{(1)}(B_R\setminus B^{(1)})=0$,

-$C^{(1)}$ is the completion of $B^{(1)}$ with respect to the geodesic distance on $C^{(1)}$.

-the volume form on $C^{(1)}$ coincides with the measure $M^{(1)}$ on $B^{(1)}$

-in a system of local coordinates,  the Riemannian metric tensor on $C^{(1)}$ reads $$g^{(1)}=\theta^{(1)}(\omega)(dx_1^2+\cdots+dx_m^2)\quad \text{ with }\quad \theta^{(1)}(\omega)=\frac{M^{(1)}(B_R)^2}{C_R^2}.$$

 %the space $(C^{(1)},g^{(1)})$ is isomorphic (as a Riemannian manifold) to the ball equipped with the Riemannian metric tensor $d^{(1)}(\omega)(dx_1^2+dx_2^2)$. That isomorphism is obtained by extending an optimal transport map $\Gamma^{(1)}:B^{(1)}\to B_R$ that pushes $\overline{M^{(1)}}$ forward to $\lambda_R$ (and we denote by $\chi^{(1)} $ its inverse).
Furthermore, from Proposition \ref{fundprop} and Theorem \ref{fundth}, $M^{(1)}\in R^{e-}_{2-\psi'(1)/n}$. This ends up the first step of the induction.

\vspace{3mm}
2) Step 2: we assume that, for some $k<n$, we may find a compact Riemannian manifold $(C^{(k)},g^{(k)})$ and a Borelian subset $B^{(k)}$ of $B_R$ such that:

-$B^{(k)}$ is dense in $B_R$ for the Euclidian distance and has full $M$-measure, that is  $M^{(k)}(B_R\setminus B^{(k)})=0$,

-$C^{(k)}$ is the completion of $B^{(k)}$ with respect to the geodesic distance on $C^{(k)}$.

-the volume form on $C^{(k)}$ coincides with the measure $M^{(k)}$ on $B^{(k)}$

-in a system of local coordinates,  the Riemannian metric tensor on $C^{(k)}$ reads $$g^{(k)}=\theta^{(k)}(\omega)(dx_1^2+\cdots+dx_m^2)\quad \text{ with }\quad \theta^{(k)}(\omega)=\frac{M^{(k)}(B_R)^2}{C_R^2}.$$

We denote by $\varphi^{(k)}=\Gamma^{(1)}\circ\cdots\circ \Gamma^{(k)}:(C^{(k)},g^{(k)})\to (B_R,d_e)$ the isometry constructed recursively. 
From Proposition \ref{fundprop}, the measure $M^{(k+1)}$ is supported by $B^{(k)}$ almost surely, that is $M^{(k+1)}(B_R\setminus 
B^{(k)})=0$ almost surely, so that $M^{(k+1)}$ extends to a measure on $C^{(k)}$ by prescribing:
$$\forall A\subset  C^{(k)} \text{ Borelian}, \quad M^{(k+1)}(A)=M^{(k+1)}(A\cap B^{(k)}).$$
Furthermore, from subsection \ref{compo}, we have $M^{(k)}\in R^{e-}_{m-\frac{k}{n}\psi'(1)}$. We now apply Theorem \ref{fundth} where the distance $d$ is equal to the geodesic distance on $C^{(k)}$, denoted by $d^{(k)}$ (the corresponding class $R_\alpha$ will be denoted by $R^{(k)}_\alpha$). Since $M^{(k)}$ is the volume form on $(C^{(k)},g^{(k)})$, we have
$$M^{(k)}(U)=Dr^m,\quad \text{for any open ball U with radius r}.$$
Hence $M^{(k)}\in R_m^{(k)-}$. The assumptions of Theorem \ref{fundth} are thus satisfies with $\alpha=m-\frac{k}{n}\psi'(1)$ and $\varsigma=m-\frac{k}{n}\psi'(1)$ (and thus we have $\psi(2)/n<\varsigma $ because of \eqref{defn}). It follows that $M^{(k+1)}\in R^{(k)-}_{m(1-\frac{\psi(2)}{mn-k\psi'(1)})+\frac{\psi(2)-\psi'(1)}{n}}$. Because of  \eqref{defn} again, we have
\begin{align*}
m(1-\frac{\psi(2)}{mn-k\psi'(1)})+\frac{\psi(2)-\psi'(1)}{n}&>m(1-\frac{\psi(2)}{mn-k\psi'(1)})\\
&>m-\frac{m\psi(2)}{n(m-\psi'(1))}\\
&>m-1,
\end{align*} so that we can show, as in the proof of Theorem \ref{gamsmallevy}, that $M^{(k+1)}$ does not charge the small sets of $(C^{(k)},g^{(k)})$. Hence we can apply Theorem \ref{th:transport} to find two optimal transport maps $\alpha^{(k+1)},\beta^{(k+1)}$ that respectively push $\overline{M}^{(k)}$ forward to $\overline{M}^{(k+1)}$ and vice versa. Furthermore,  $\beta^{(k+1)}$ can be rewritten as $(\varphi^{(k)})^{-1}\circ \Gamma^{(k+1)} \circ \varphi^{(k)}$ where $\Gamma^{(k+1)}:B_R\to B_R$ is the gradient of some convex function. The function $ \varphi^{(k+1)}=\varphi^{(k)}\circ \beta^{(k+1)}=\Gamma^{(k+1)}\circ\cdots\circ \Gamma^{(1)}$ thus pushes $M^{(k+1)}$ forward to $\lambda_R$. Besides the quantity 
 $$\inf_{\substack{T:B_R\to B_R\\T_{\#}\overline{M}^{(k+1)}=\overline{M}^{(k)}}} \int_{B_R}d^{(k)}(T(x),x)^2\overline{M}^{(k+1)}(dx)$$ is achieved  at $\beta^{(k+1)}$. Since $d^{(k)}(x,y)=d_e(\varphi^{(k)}(x),\varphi^{(k)}(y))$, we deduce that $\beta^{(k+1)}$ realizes the above infimum if and only if $ \varphi^{(k+1)}=\varphi^{(k)}\circ \beta^{(k+1)} $ realizes the infimum
  $$\inf_{\substack{T:B_R\to B_R\\T_{\#}\overline{M}^{(k+1)}=\lambda_R}} \int_{B_R}|T(x)-\varphi^{(k)}(x)|^2\overline{M}^{(k+1)}(dx).$$
Then we can apply the same machinery as in subsection \ref{optimal} to construct the $k+1$-th Riemannian structure, which ends up the induction.\qed

 \appendix
 
 % 
%%%%%%%%%%%%%%%%%%%%%%%%%%%%%%%%%%%%%%%%%%%%%%%%%%%%%%%%%%%%%%%%%%%%%%%%%%%%

\section{Background about optimal transport theory}\label{optimaltheory}

%%%%%%%%%%%%%%%%%%%%%%%%%%%%%%%%%%%%%%%%%%%%%%%%%%%%%%%%%%%%%%%%%%%%%%%%%%%
\subsection{Monge problem}
%%%%%%%%%%%%%%%%%%%%%%%%%%%%

We remind the reader of  the following classical results
\begin{definition}{\bf Push-forward of measures.}
Let $\mu,\nu$ be two measures respectively defined on the measured spaces $E$ and $F$. We will say that a measurable mapping $\varphi:F\to E$ pushes the measure $\nu$ forward to $\mu$ if both measures $\mu$ and $\nu\circ \varphi^{-1}$ coincide. In that case, we write $\varphi_\#\nu=\mu$.
\end{definition}
\begin{definition}{\bf Small sets.}
Given a metric space $(X,d)$ with Hausdorff dimension $n$, a small set is a set with Hausdorff dimension not greater than $n-1$.
\end{definition}
%\begin{definition}{\bf $d^2/2$ convex functions.}
%A function $\psi:\T^d\to \R\cup\{+\infty\}$ is said to be $d^2/2$-convex if it is not identically $+\infty$ and if there exists $\zeta:\T^d\to \R\cup\{\pm \infty\}$ such that
%$$\forall x\in \T^d,\quad \psi(x)=\sup_{y\in\T^d}\big(\zeta(y)-d^2(x,y)/2\big).$$
%\end{definition}

 Given two probability measures $\mu$ and $\nu$ on $B_R$, a coupling of $(\mu,\nu)$ is a probability measure $\pi$ on $B_R\times B_R$ with marginals $\mu$ and $\nu$. A coupling $\pi $ is said to be deterministic if there is a measurable map $T:B_R\to B_R$ such that the map $x\in B_R\mapsto (x,T(x))$ pushes $\mu$ forward to $\pi$. In particular, for all $\nu$-integrable function $\varphi$, one 
 has
 $$\int_{B_R}\varphi(y)\,d\nu(y)=\int_{B_R}\varphi(T(x))\,d\mu(x).$$ Such a map $T$ is called a transport map between $\mu$ and $ \nu$.
 
 The Monge-Kantorovich problem on the ball $B_R$ can be formulated as follows. Given a cost function $c$ defined on $B_R\times B_R$, one looks for a coupling $\pi $ of $(\mu,\nu)$ that realizes the infimum
$$C(\mu,\nu)=\inf\int_{B_R\times B_R}c(x,y)d\pi(x,y) $$ where the infimum runs over all the coupling $\pi$ of $(\mu,\nu)$. Such a coupling is called optimal transference plan. If the coupling $\pi$ is deterministic, the corresponding transport map $T$ is called optimal transport map. The optimal transport cost is then the value
$$\int_{B_R\times B_R}c(x,T(x))d\mu(x) .$$
The search of deterministic optimal transference plans is called the Monge problem.  

%\subsection{$c$-convexity}
%%%%%%%%%%%%%%%%%%%%%%%%%%%%%
%Let $(\mathcal{X},d)$  be a Polish space and $c:(x,y)\in \mathcal{X}\times \mathcal{X}\mapsto d(x,y)^2$.  A function $\psi:\mathcal{X}\to (-\infty;+\infty]$ is said to be $c$-convex if it is not identically $+\infty$ and if there is a function $\phi:\mathcal{X}\to  [-\infty;+\infty]$ such that 
%$$\forall x\in  \mathcal{X},\quad \psi(x)=\sup_{y\in \mathcal{X}}\big(\phi(y)-d(x,y)^2\big).$$
%Then its $c$-transform is the function $\psi^c$ defined by 
%$$\forall y\in  \mathcal{X},\quad \psi^c(y)=\inf_{x\in \mathcal{X}}\big(\psi(x)+d(x,y)^2\big),$$
%and its $c$-subdifferential at the point $x\in  \mathcal{X}$ is given by
%$$\partial_c\psi(x)=\{y\in  \mathcal{X};\,\,\psi^c(y)-\psi(x)=d(x,y)^2.$$

\subsection{Solution to the Monge problem}
%%%%%%%%%%%%%%%%%%%%%%%%%%%%

We have (see \cite[Theorem 10.28]{villani} and \cite[Theorem 2.12 iv]{villani2} for the last statement)
\begin{theorem}\label{th:transport}
Let $\mathcal{X}$  be a Riemannian manifold isometric (as a smooth Riemannian manifold) to the closed ball $B_R$. We denote by $f: \mathcal{X}\to B_R$  the corresponding isometry of Riemannian structures.
Let $c:\mathcal{X}\times \mathcal{X}\to \R$ be the cost function given by 
$$c(x,y)=d(x,y)^2 $$ and $\mu,\nu$ two probability measures on $\mathcal{X}$. Assume that the measure $\mu$ does not give mass to small sets. Then:

1) there is a unique (in law) optimal coupling $\pi$ of $(\mu,\nu)$ and it is deterministic. 

2) there is a unique optimal transport map $T$ (i.e. uniquely determined $\mu$ almost everywhere) solving the Monge problem. Furthermore, we can find a lower semi-continuous convex function $\phi$ defined on $B_R$ such that $$T(x)=f^{-1}\circ \nabla\phi \circ f(x)$$ for every $x\in f^{-1}\big(\{y\in \R;\phi\,\,\text{ is differentiable at }y\}\big).$

3) ${\rm Supp}(\nu)=\overline{T({\rm Supp}(\mu))}$.

4) Finally, if $\nu$ does not give mass to small sets either, then there is also a unique optimal transport map $T'$ solving the Monge problem
(of pushing $\nu$ forward to $\mu$). We can also find a lower semi-continuous convex function $\psi$ defined on $B_R$ such that $$T'(x)=f^{-1}\circ \nabla\phi \circ f(x)$$ for every $x\in f^{-1}\big(\{y\in \R;\psi\,\,\text{ is differentiable at }y\}\big).$
$T$ and $T'$ satisfy, for $\mu$ almost every $x\in\mathcal{X}$ and $\nu$ almost every $y\in\mathcal{X}$,
$$T'\circ T(x)=x,\quad  T\circ T'(y)=y.$$
\end{theorem}

\noindent {\it Proof.} There is an easy way to deduce the above theorem from \cite[Theorem 10.28]{villani}. Because of the isometry with the closed ball $B_R$, the above theorem is basically of Euclidian nature. Indeed, it is plain to see that $\pi$ is an optimal coupling  of $\mu,\nu$ for the cost function $c(x,y)=d(x,y)^2$ on $\mathcal{X}\times \mathcal{X}$ if and only if $\Pi=\pi\# (f,f)$ is a coupling of the probability measures of $f\#\mu,f\#\nu$ on $B_R$. In the same way, $T:\mathcal{X}\to \mathcal{X}$ is an optimal transport map such that $T\#\mu=\nu$ if and only if $\theta=f\circ T \circ f^{-1}$ is an optimal transport map such that $\theta\#(f\#\mu)=f\#\nu$ for the Euclidian quadratic cost on $B_R$. So the proof of the above theorem boils down to applying \cite[Theorem 10.28]{villani} in the Euclidian case with quadratic cost function. It is then plain to complete the proof.\qed

\begin{remark}
In  case 4) is satisified, it is more convenient to restrict the support of $T$ and $T'$ respectively to $\{x\in \mathcal{X}; \phi\,\,\text{ is differentiable at }f(x)\text{ and } T'\circ T(x)=x \}$ and  $\{x\in \mathcal{X}; \psi\,\,\text{ is differentiable at }f(x)\text{ and } T\circ T'(x)=x \}$. In that way, $T:{\rm Supp}(T)\to{\rm Supp}(T')$ and $T':{\rm Supp}(T')\to{\rm Supp}(T)$ are both bijections.
\end{remark}

\subsection*{Acknowledgements}
The authors are thankful to G. Carlier, B. Nazaret, F. Santambrogio and C. Villani for numerous fruitful discussions about optimal transportation theory.

%%%%%%%%%%%%%%%%%%%%%%%%%%%%%%%%%%%%%%%%%%%%%%%%%%%%%%%%%%%%%%%%%%%%%%%%%%%%%%%%%%%%%%%%%%

\end{document}